\documentstyle[12pt,amscd,amssymb]{amsart}
\newtheorem{pr}{Proposition}
\newtheorem{lm}{Lemma}[subsection]

\newcommand{\proj}{{\Bbb P}}
\newcommand{\grass}{{\Bbb G}}

\newcommand{\com}{{\Bbb C}}
\newcommand{\LL}{{\cal{L}}}
\newcommand{\NN}{{\cal{N}}}
\newcommand{\TT}{{\cal{T}}}
\newcommand{\ZZ}{{\cal{Z}}}
\newcommand{\barr}{\overline}
\newcommand{\rarr}{\rightarrow}
\newcommand{\oh}{{\cal{O}}}
\newcommand{\Q}{{\Bbb{Q}}}
\newcommand{\deli}{\bigtriangleup}
\newcommand{\eqq}{\stackrel {\sim}{=}}

\newcommand{\M}{\barr{M}}
\newcommand{\N}{\barr{N}}

\newcommand{\HH}{{\cal{H}}}
\newcommand{\tl}{\tilde}

\begin{document}
\title{Intersections of $\Bbb{Q}$-Divisors on Kontsevich's
Moduli Space $\M_{0,n}(\proj^r,d)$ and Enumerative Geometry}
\maketitle
\begin{center}
6 April 1995
\end{center}
\begin{center}
Rahul Pandharipande
\end{center}
\pagestyle{plain}
\baselineskip=14pt
\setcounter{section}{-1}
\section{{\bf Introduction}}
Let $\com$ be the field of complex numbers.
Let $(C, p_1, \ldots, p_n)$ be a connected, reduced, projective, nodal curve
over
$\com$ with  $n$ nonsingular marked points $(p_1, \ldots, p_n)$.
Let $\omega_C$ be the dualizing sheaf of $C$.
An algebraic map $\mu: (C, p_1,\ldots, p_n) \rarr \proj^r$ is
{\em Kontsevich stable} if $\omega_C(p_1+ \ldots +p_n) \otimes
\mu^*(\oh_{\proj^r}(3))$ is ample on $C$.
Let $\M_{g,n}(r,d)$ be the coarse moduli space of degree $d$, Kontsevich
stable maps from $n$-pointed, genus $g$ curves to $\proj^r$.
In the genus zero case, $\M_{0,n}(r,d)$ is an irreducible,
projective variety with finite quotient singularities.
Only the following cases will be considered here:
$$d\geq 0,\ \ g=0,\ \ r\geq 2.$$
The stack of Kontsevich stable maps was first defined in
[K-M] and [K]. A treatment of the corresponding coarse moduli spaces
can also be found in [P] and [Al].

The dimension of $\M=\M_{0,n}(r,d)$ is $m=rd+d+r+n-3$.
Let $Pic(\M)$ be the Picard group of
line bundles.
Let $A_{m-1}(\M)$ be the Chow group of Weil divisors modulo
rational equivalence.
Since $\M$ has finite quotient singularities,
every Weil divisor is $\Bbb{Q}$-Cartier. Therefore, there
is a canonical isomorphism:
$$Pic(\M)\otimes \Q \rarr A_{m-1}(\M) \otimes \Q\ .$$
$Pic(\M)\otimes \Q$ is a finite dimensional vector space.
An explicit set of generators is given below.

Let $P=\{1,2,\ldots n\}$ be the set of markings ($P$ may be the empty set).
The $n$ markings of the moduli problem yield $n$ canonical line bundles
$\LL_i= \nu_i^*(\oh_{\proj^r}(1))$ on
$\M$ via the $n$ evaluation maps
$\forall i \in P, \ \ \nu_i: \M \rarr \proj^r$.
The {\em boundary} of $\M$ is the locus corresponding to
maps with reducible domain curves. Since the  boundary is
of pure codimension $1$ in $\M$,
each irreducible component is a Weil divisor. The irreducible components
of the boundary are in bijective correspondence with
data of weighted partitions $(A\cup B, d_A, d_B)$ where:
\begin{enumerate}
\item[(i.)] $A\cup B$ is a partition of $P$.
\item[(ii.)] $d_A+d_B=d$, $d_A>0$, $d_B>0$.
\item[(iii.)] If $d_A=0$ (resp. $d_B=0$), then $|A|\geq 2$ (resp.
$|B| \geq 2$).
\end{enumerate}
For example, if $P=\emptyset$, then $A=B=\emptyset$ and
the boundary components
correspond to positive partitions $d_A+d_B=d$. Let $\deli$ be the
set of components of the boundary.

In case $d\geq 1$,
a Weil divisor
is obtained on $\M$ by considering the locus of
$\M$ corresponding to maps meeting a fixed $r-2$ dimensional
linear subspace of $\proj^r$ (note $r\geq 2$). It is shown in
[P] that this incidence Weil divisor is actually Cartier. Denote
the corresponding line bundle on $\M$ by $\HH$. For convenience,
let $\HH=0$ in case $d=0$.

\begin{pr} Results on generation:
\begin{enumerate}
\item[(i.)] If $d=0$, $g=0$, $r\geq 2$,
$\{\LL_i\}\cup \deli$ generate $Pic(\M)$.
\item[(ii.)]If $d\geq 1$, $g=0$, $r\geq 2$,
$\{\LL_i\} \cup \deli \cup \{\HH\}$
generate $Pic(\M)\otimes \Q$.
\end{enumerate}
\label{gen}
\end{pr}
\noindent
If $d=0$, then (by stability) $n\geq 3$
 and $\M_{0,n}(r,0) \eqq \barr{M}_{0,n}\times \proj^r$
where $\barr{M}_{0,n}$ is the Mumford-Knudsen space.
In this case, $\LL_i$ is the pull-back of $\oh_{\proj^r}(1)$ from the
second factor. Therefore, part (i) is a consequence of the boundary generation
of
$Pic(\barr{M}_{0,n})$.

There is an intersection pairing $A_1(\M) \otimes Pic(\M)
 \rarr \Bbb{Z}$.
Let $Null \subset Pic(\M)$ be the null space with respect to the
intersection pairing. Define
$$Num(\M) = Pic(\M)/Null.$$
By Proposition (\ref{gen}), the classes $\{\LL_i\} \cup
\deli \cup \{\HH\}$ generate $Num(\M)\otimes \Q$.
The relations between these generators in $Num(\M)\otimes \Q$
can be algorithmically determined by calculating intersections
with curves. It will be shown that all the relations in
$Num(\M)\otimes \Q$ are obtained from linear equivalences in
$Pic(\M)\otimes \Q$.
\begin{pr}
\label{dimm}
The canonical map $Pic(\M)\otimes \Q \rarr Num(\M)\otimes \Q$ is
an isomorphism. The Picard numbers are:
\begin{enumerate}
\item[{}]$(n=0)$, $\ dim_{\Q}\ Pic(\M)\otimes \Q =  [{d\over 2}]+1.$
\item[{}]$(n\geq 1)$,  $\ dim_{\Q} \ Pic(\M) \otimes \Q=
(d+1)\cdot 2^{n-1}-{n\choose 2}.$
\end{enumerate}
\end{pr}
The main result of this paper concerns the computations of
top intersection products in $Pic(\M)\otimes \Q$.
\begin{pr}
\label{top}
Let $d\geq 0$, $g=0$, $r\geq 2$.
There exists an explicit algorithm for calculating
the top dimensional intersection products of the $\Q$-Cartier divisors
$\{\LL_i\} \cup \deli \cup \{\HH\}$ on $\M$.
\end{pr}
Consider the space $R(d,r)$ of degree $d\geq 1$ rational curves in $\proj^r$
($r\geq 2$). The dimension of $R(d,r)$ is $rd+r+d-3$.
Classically, the characteristic numbers of $R(d,r)$ are
the numbers of degree $d$ rational curves in $\proj^r$
passing through $\alpha_i$ general linear spaces of
codimension $i$ (for $2\leq i \leq r$) and tangent to
$\beta$ general hyperplanes where
$$(i-1)\cdot \alpha_i + \beta = dim \ R(d,r).$$
The characteristic numbers of rational curves excluding
tangencies ($\beta=0$) have been determined recursively  by M. Kontsevich
and Y. Manin in [K-M] (also by Y. Ruan and G. Tian in [R-T]).
The divisor in $\M$ corresponding to
the hyperplane tangency condition can be expressed as a linear
combination of the classes $\{\LL_i\} \cup \deli \cup \{\HH\}$.
The characteristic
numbers can then be expressed as top intersection products
of  $\{\LL_i\} \cup \deli \cup \{\HH\}$ on suitably
chosen Kontsevich spaces of maps $\M$.
Therefore, all the characteristic numbers can be
calculated by Proposition (\ref{top}).
\begin{pr}
\label{charny}
There exists an explicit algorithm for calculating
all the characteristic numbers of rational curves in
projective space.
\end{pr}
P. Di Francesco and C. Itzykson have modified the methods of [K-M]
to determine some ($\beta\neq 0$) characteristic numbers
for rational plane curves ([D-I]). Unfortunately, the relations they obtain
from the WDVV associativity equations do not suffice to
recursively determine all the characteristic numbers for rational
plane curves from a finite set of data.

The structure of the paper is as follows. Propositions (\ref{gen}) and
(\ref{dimm}) are proven in section (\ref{cone}). In section (\ref{calc}),
several geometric classes are explicitly computed in $Pic(\M)\otimes
\Q$. These classes will be used in the algorithms of Propositions
(\ref{top}) and (\ref{charny}). The algorithms are established in section
(\ref{inter}). Section (\ref{exam}) is devoted to calculations of some
characteristic numbers of rational curves for small values of
$(d,r)$. As a final application, a new formula for cuspidal
rational curves is derived in section (\ref{cusp}).

The problem of calculating tangency characteristic numbers
via Kontsevich's moduli space was suggested to the author by W. Fulton.
Conversations with W. Fulton on related topics have been of
significant aid. Thanks are due to S. Kleiman and P. Aluffi for
mathematical and historical remarks.
The remarkable ideas in [K-M] have been a
source of inspiration.

\section{$Pic(\M)\otimes \Q$ And $Num(\M)\otimes \Q$}
\label{cone}
\setcounter{subsection}{-1}
\subsection{Summary}
Propositions (\ref{gen}) and (\ref{dimm}) are established in
sections (\ref{jenny}) and (\ref{realy}) respectively. Since these
results are well known for $d=0$,
$$\M_{0,n}(r,0) \eqq \barr{M}_{0,n}\times
\proj^r,$$
the conditions $d\geq 1$, $g=0$, $r\geq2$ are assumed throughout
sections (\ref{jenny}) and (\ref{realy}).
\subsection{Generators}
\label{jenny}
The proof of Proposition (\ref{gen}) is divided into four cases
depending upon the number $n$ of marked point.
\begin{lm} If $n\geq 3$, then $Pic(\M)\otimes
\Q$ is generated by $\deli \cup \{\HH\}.$
\label{n3}
\end{lm}
\begin{pf}
Let $V=\bigoplus_{0}^{r} H^0(\proj^1, \oh_{\proj^1}(d))$.
Let $U\subset \proj(V)$ be the Zariski open set corresponding to
a well defined (basepoint free) degree $d$ map from $\proj^1$ to $\proj^r$.
The complement of $U$ in $\proj(V)$ is of codimension at least $r\geq 2$.
There is a universal map
$$\proj^1 \times U \rarr \proj^r.$$
Fix the first three marked points to be $0$, $1$, $\infty\in \proj^1$.
Let
$$W= \proj^1 \times \ldots \times \proj^1 \ \setminus \{D_{i,j}, \ S_{0,i},
\ S_{1,i}, \ S_{\infty, i} \}$$
where the product is taken over $n-3$ factors.
$D_{i,j}$ is the large diagonal determined by factors
$i$ and $j$.
$S_{0,i}$ is the locus where the $i^{th}$ factor is $0\in \proj^1$.
$S_{1,i}$, $S_{\infty, i}$ are defined similarly.
It follows there is a universal family of Kontsevich stable
degree $d$ maps of $n$-pointed curves:
$$\proj^1 \times W \times U \rarr \proj^r.$$
The maps of the family are
automorphism-free and distinct.
By the universal property, there is an injection
$W\times U \rarr \M$. A tangent space calculation
shows $W\times U$ is an open set of $\M$. The complement of
$W\times U$ is the boundary of $\M$.
Hence $A_{m-1}(\M)$ is generated by $\deli$ and $A_{m-1}(W\times U)$.
Information about $A_{m-1}(W\times U)$ is obtained from
the open inclusion
\begin{equation}
\label{inclu}
W\times U \subset \proj^1 \times \ldots \times \proj^1 \times \proj(V).
\end{equation}
The Picard group of the right side of (\ref{inclu}) is generated
by the pull-backs of $\oh(1)$ from each factor. The pull-backs
from the $\proj^1$ factors are trivial on $W\times U$ because
of the removal of the loci $S_{0,i}$. Hence,
$A_{m-1}(W\times U)$ is generated by $\oh_{\proj(V)}(1)$.
It is easily seen $\HH$ restricted to $W\times U$
is the pull-back of a resultant hypersurface in $\proj(V)$.
Therefore, $\HH$ restricted to $W\times U$
is
linearly equivalent to a multiple of $\oh_{\proj(V)}(1)$.
\end{pf}

There are canonical morphisms $\M_{0,n}(r,d) \rarr \M_{0,n-1}(r,d)$
obtained by omitting the last marked point. Results for $0\leq n \leq 2$
are obtained via these morphisms.
\begin{lm} If $n=2$, then $Pic(\M)\otimes \Q$ is generated by
$\deli \cup \{\LL_1, \LL_2\}$.
\end{lm}
\begin{pf}
Let $\N=\M_{0,3}(r,d)$ and $\M=\M_{0,2}(r,d)$. Fix a
hyperplane $H_3\subset \proj^r$. Let $X= \nu_3^{-1}(H_3)$
where $\nu_3$ is the third evaluation map,
$\nu_3: \N \rarr \proj^r.$
There is a map $\rho: X \rarr \M$ obtained by omitting the
third point. The map $\rho$ is surjective and generically
finite. Let $Z\subset \M$ be the open set
corresponding to Kontsevich stable maps satisfying the
following conditions:
\begin{enumerate}
\item[(i.)] The domain curve is $\proj^1$.
\item[(ii.)] The images of the marked points $\{1,2\}$ do not lie in $H_3$.
\end{enumerate}
It is clear the the complement of $Z$ is the boundary union
$\nu_1^{-1}(H_3)$, $\nu_2^{-1}(H_3)$.
By the definition of $Z$, $\rho^{-1}(Z) \rarr Z$ is
a finite, projective morphism.
If $A_{m-1}(\rho^{-1}(Z))=0$, then $A_{m-1}(Z)$ is torsion.
To establish the Lemma,
it therefore suffices to prove $A_{m-1}(\rho^{-1}(Z))=0$.

In the notation of the proof of Lemma (\ref{n3}),
$\rho^{-1}(Z) \subset U\subset \M_{0,3}(r,d)$. In fact, the following
is easily seen:
$$\rho^{-1}(Z)  = \ U \cap L_{\infty}(H_3) \ \setminus \ \{L_0(H_3),
L_1(H_3)\}.$$
$L_p(H_3)$ is the hyperplane in $U$ corresponding to maps sending the
point $p\in \proj^1$ to $H_3$.
$U\cap L_{\infty}(H_3)$ is a an open set of $L_{\infty}(H_3)$
with complement of codimension at least $2$. Hence,
$A_{m-1}(U\cap L_{\infty}(H_3))=\Bbb{Z}$ generated by the hyperplane
class.
Since $\rho^{-1}(Z)\subset U\cap L_{\infty}(H_3)$ is the
complement of hyperplanes, the desired conclusion
$A_{m-1}(\rho^{-1}(Z))=0$ is obtained.
\end{pf}

\begin{lm} If $n=1$, then $Pic(\M)\otimes \Q$ is generated by
$\deli \cup \{\LL_1, \HH\}$.
\end{lm}
\begin{pf}
Let $\N=\M_{0,3}(r,d)$ and $\M=\M_{0,1}(r,d)$. Fix two hyperplanes
$H_2, H_3\subset \proj^r$. Let $X= \nu_2^{-1}(H_2)\cap
\nu_3^{-1}(H_3)$
where $\nu_2$, $\nu_3$ are  the second and third evaluation maps on $\N$.
There is a map $\rho: X \rarr \M$ obtained by omitting the second and
third points. The map $\rho$ is surjective and generically
finite. Let $Z\subset \M$ be the open set
corresponding to Kontsevich stable maps satisfying the
following conditions:
\begin{enumerate}
\item[(i.)] The domain curve is $\proj^1$.
\item[(ii.)] The image of the marked point $\{1\}$ does not lie in $H_2 \cup
H_3$.
\item[(iii.)] The map does not pass through the intersection
$H_2 \cap H_3$.
\end{enumerate}
The complement of $Z$ is the boundary union
$\nu_1^{-1}(H_2)$, $\nu_1^{-1}(H_3)$, and $D_{2,3}$. $D_{2,3}$
is the Cartier divisor of maps passing through $H_2 \cap H_3$.
$D_{2,3}$ is a divisor in the linear series of $\HH$.
By the definition of $Z$, $\rho^{-1}(Z) \rarr Z$ is
a finite, projective morphism.
As before,
it suffices to prove $A_{m-1}(\rho^{-1}(Z))=0$.

Let $S\subset U$ be the union of the hyperplane sections
$\{L_0(H_2), L_0(H_3)\}$ with the resultant hypersurface
of maps meeting $H_2\cap H_3$.
Conditions (i), (ii), and (iii) imply:
$$\rho^{-1}(Z)  = \ U \cap L_{1}(H_2) \cap L_{\infty}(H_3) \
\setminus  \ S.$$
As before, $A_{m-1}(U\cap L_1(H_2) \cap L_{\infty}(H_3))=\Bbb{Z}$
generated by the hyperplane class. $S$ is a union
of hyperplane classes and multiples of hyperplane classes. Hence,
$A_{m-1}(\rho^{-1}(Z))=0$.
\end{pf}

\begin{lm} If $n=0$, then $Pic(\M)\otimes \Q$ is generated by
$\deli \cup \{\HH\}$.
\label{n0}
\end{lm}
\begin{pf}
Let $\N=\M_{0,3}(r,d)$ and $\M=\M_{0,0}(r,d)$. Fix three general hyperplanes
$H_1, H_2, H_3\subset \proj^r$. Let $X=  \nu_1^{-1}(H_1)
\cap \nu_2^{-1}(H_2)\cap
\nu_3^{-1}(H_3)$
where the $\nu_i$ are evaluation maps on $\N$.
There is a map $\rho: X \rarr \M$ obtained by omitting the marked points.
The map $\rho$ is surjective and generically
finite. Let $Z\subset \M$ be the open set
corresponding to Kontsevich stable maps satisfying the
following conditions:
\begin{enumerate}
\item[(i.)] The domain curve is $\proj^1$.
\item[(ii.)] The map does not pass through the intersections
$H_1\cap H_2$, $H_1\cap H_3$, or $H_2 \cap H_3$.
\end{enumerate}
The complement of $Z$ is the boundary union
$D_{1,2}$, $D_{1,3}$, and $D_{2,3}$.
By the definition of $Z$, $\rho^{-1}(Z) \rarr Z$ is
a finite, projective morphism.
As before,
it suffices to prove $A_{m-1}(\rho^{-1}(Z))=0$.

Let $S\subset U$ be the union of the three resultant hypersurfaces
of maps meeting $H_1\cap H_2$, $H_1\cap H_3$, and $H_2\cap H_3$.
Let $I\subset U$ be the hyperplane intersection
defined by $I=U\cap L_0(H_1)\cap L_1(H_2) \cap L_{\infty}(H_3)$.
Conditions (i) and (ii) imply:
$$\rho^{-1}(Z) = \ I
\ \setminus \ S \cap I.$$
Note $S\cap I$ contains the
intersections of the following hyperlanes with $I$:
$$\{L_0(H_2), L_0(H_3), L_1(H_1), L_1(H_3), L_{\infty}(H_1)
, L_{\infty}(H_2)\}.$$
As before, $A_{m-1}(U\cap I)=\Bbb{Z}$
generated by the hyperplane class. Since $S\cap I$ is a union of
of hyperplane classes and multiples of hyperplane classes,
$A_{m-1}(\rho^{-1}(Z))=0$.
\end{pf}
Lemmas (\ref{n3}) - (\ref{n0}) yield Proposition (\ref{gen}).

\subsection{Relations}
\label{realy}
Curves in $\M=\M_{0,n}(r,d)$ are easily found.
The following construction will be required for the calculations
below.
Let $C$ be a nonsingular, projective curve.
Let $\pi: S=\proj^1 \times C \rarr C$. Select $n$ sections
$s_1, \ldots, s_n$
of $\pi$.
A point $x\in S$ is an
{\em intersection point} if two or more sections contain $x$.
Let $\NN$ be a line bundle on $S$ of type $(d,k)$ where
$k$ is very large. Let $z_l\in H^0(S, \NN)$  $(0\leq l \leq r)$
determine a rational map $\mu: S - \rarr \proj^r$ with simple
base points. A point $y\in S$ is a {\em simple base point} of
degree $1\leq e\leq d$ if  the blow-up of $S$ at $y$ resolves $\mu$ locally
at $y$ and the resulting map is of degree $e$ on the exceptional
divisor $E_y$.
The set of {\em special points} of $S$ is the
union of the intersection points and the simple base points.
Three conditions are required:
\begin{enumerate}
\item [(1.)] There is at most one special point in each fiber of
$\pi$.
\item [(2.)] The sections through each intersection point $x$
have distinct tangent directions at $x$.
\item [(3.)] If $n$ or $n-1$ sections pass through the
point $x\in S$, then $x$ is not a simple base point of degree $d$.
(If $n=0$ or $1$, there are no simple base points of degree $d$.)
\end{enumerate}
Let $\barr{S}$ be the blow-up of $S$ at the special points.
It is easily seen $\barr{\mu} : \barr{S} \rarr \proj^r$ is
Kontsevich stable family of $n$-pointed, genus $0$ curves over
$C$. Condition (2) ensures the strict transforms of the sections
are disjoint. Condition (3) implies Kontsevich stability.
There is a canonical morphism $C \rarr \M$.
Condition (1) implies $C$ intersects the boundary components transversally.

\begin{lm}
\label{ihh}
Results on the span of $\{\HH, \LL_1, \LL_2\}$:
\begin{enumerate}
\item[(i.)] The element $\HH$ is not contained in the
linear span of $\deli$ in $Pic(\M)\otimes \Q$.
\item[(ii.)] If $n=1$, $\{\HH, \LL_1\}$ are independent modulo $\deli$.
\item[(iii.)] If $n=2$, $\{\LL_1, \LL_2\}$ are independent modulo $\deli$.
\end{enumerate}
\end{lm}
\begin{pf}
Consider $\pi: S=\proj^1 \times C \rarr C$ with $n$ trivial sections.
There are no intersection points.
Let $\NN$, $z_l\in H^0(S, \NN)$ be such that $\mu$ has no
base points (note: since $r\geq 2$, this is easily accomplished).
$\NN$ has degree type $(d,k)$.
For each
component $K\in \deli$, $C\cdot K=0$.
A simple calculation yields $C\cdot \HH= \NN\cdot \NN=2dk$.
Hence $\HH$ is not contained in the span of $\deli$.

Consider  $\pi: S=\proj^1 \times \proj^1 \rarr
\proj^1$. Let $s$ be the trivial section; let $s'$ be the
diagonal section. Let $\mu: S \rarr \proj^r$ be a base point
free map of type $(d,k)$. The two sections $s$, $s'$ determine
two maps $\tau, \tau': \proj^1 \rarr \M_{0,1}(r,d)$.
Intersection via $\tau$ yields:
$$\proj^1 \cdot \HH= 2dk, \ \ \proj^1 \cdot \LL_1= k.$$
Intersection via $\tau'$ yields:
$$\proj^1 \cdot \HH=2dk, \ \ \proj^1 \cdot \LL_1= d+k.$$
In both cases $\proj^1\cdot K=0$ for any $K\in \deli$.
Therefore $\{\HH, \LL_1\}$ are independent modulo $\deli$
in $Pic(\M) \otimes \Q$ for $n=1$.

In the $n=2$ case, twisted families must be considered.
Let $E(a,b)$ be the rank two bundle
$\oh_{\proj^1}(a)\oplus \oh_{\proj^1}(b)$ over $\proj^1$.
Let $S(a,b)= \proj(E(a,b))$. Let $$\NN= \oh_{\proj(E)}(d)\otimes
\pi^*(\oh_{\proj^1}(k)).$$
For large $k$, let $\mu: S(a,b)\rarr \proj^r$ be a base point
free map. The sub-bundles $\oh(a)$, $\oh(b)$ define
sections $s_1$ and $s_2$. There is an induced map
$\proj^1 \rarr \M_{0,2}(r,d)$.
A calculation yields:
$$ \proj^1\cdot \LL_1= -ad+k, \ \ \proj^1 \cdot \LL_2=-bd+k.$$
As before $\proj^1\cdot K=0$ for any $K\in \deli$.
It follows $\{\LL_1, \LL_2\}$ are independent modulo
$\deli$ in $Pic(\M) \otimes \Q$ for $n=2$.
\end{pf}

If  $n\geq 1$,
let $\deli_i\subset \deli$ be the subset of boundary components
$(A\cup B, d_A, d_B)$ with marking  partition $|A|+|B|=n$
equal to the partition $i+(n-i)=n$. There is a disjoint union
$$\deli = \bigcup_{i=0}^{[{n\over 2}]} \deli_i.$$
Let $\deli'=\deli \setminus( \deli_0 \cup \deli_1)$.
\begin{lm}
\label{i01} Results on the span of $\deli_0$, $\deli_1$:
\begin{enumerate}
\item[(i.)]
If $n=0$, $\deli_0=\deli$ is a set of linearly independent elements of
$Pic(\M)\otimes \Q$.
\item[(ii.)] If $n=1$, $\deli_0=\deli_1=\deli$ is a set
of linearly independent elements of $Pic(\M)\otimes \Q$.
\item[(iii.)] If $n\geq 2$,
$\deli_0 \cup \deli_1$ is a set of linearly independent elements
of $Pic(\M)\otimes \Q$. Moreover, the span of $\deli_0 \cup \deli_1$
does not intersect the span of $\deli'$ in $Pic(\M)\otimes \Q$.
\end{enumerate}
\end{lm}
\begin{pf}
Let $\pi: S=\proj^1\times C \rarr C$ be as above with $n$
trivial sections. Let $\NN$ be a line bundle on $S$ of degree
type $(d,k)$.
For large degrees $k$, the simple
base points of $\mu$ of degree $1\leq e \leq d$
can be selected arbitrarily satisfying
conditions (1) and (3).
For suitable choices of simple base points and base point degrees on $S$,
the classes in assertions (i-iii) can be seen to be independent
in $Num(\M) \otimes \Q$. Therefore, the classes are independent in
$Pic(\M)\otimes \Q$.
\end{pf}

In case $n=0$, $\deli \cup \{\HH\}$ is a basis of $Pic(\M)\otimes Q$
via Lemmas (\ref{n0}), (\ref{ihh}), and (\ref{i01}).
For $1\leq n \leq 3$, $\deli_0 \cup \deli_1= \deli$.
Hence, the Lemmas show the generators of section (\ref{jenny})
are also bases for $1\leq n \leq 3$. The Picard numbers of
Proposition (\ref{dimm}) can be verified for $0\leq n \leq 3$.

For $n\geq 4$, let $\barr{M}_{0,n}$ be the Mumford-Knudsen
moduli space of $n$-pointed, genus $0$ curves. The boundary
components of $\barr{M}_{0,n}$ correspond bijectively to
partitions $A\cup B$ of $P=\{1,2,\ldots, n\}$ such that
$|A|, |B|\geq 2$. The boundary components generate
$Pic(\barr{M}_{0,n})$. The three boundary components of
$\barr{M}_{0,4}$ are linearly equivalent.
A four element subset $Q\subset P$
induces a natural map $\barr{M}_{0,n} \rarr \barr{M}_{0,Q}$.
The pull-backs of the basic boundary linear equivalences on
$\barr{M}_{0,Q}$ induces boundary linear equivalences
on $\barr{M}_{0,n}$.
The relations among the boundary
components of $\barr{M}_{0,n}$
are generated by these pull-back linear equivalences as $Q$ varies among all
four element subsets of $P$.
$Pic(\barr{M}_{0,n})$ is a free group of rank
$$2^{n-1}-{n-1 \choose 2}-n.$$ Since there are
$${2^{n}-2-2n \over 2}= 2^{n-1}-1-n$$
boundary components of $\barr{M}_{0,n}$, it follows there are
${n-1\choose 2}-1$ independent relations among the boundary
components. Finally, $Pic(\barr{M}_{0,n}) \eqq Num(\barr{M}_{0,n})$.
See [Ke] for proofs of these
results.

Let $n\geq 4$.
There is canonical morphism $\eta:\M=\M_{0,n}(r,d)\rarr \barr{M}_{0,n}$.
The $\eta$ pull-back of a boundary component of $\barr{M}_{0,n}$
is a non-empty, multiplicity-free sum of boundary components $\deli'$ of $M$:
$$\eta^{-1}\big( (A\cup B) \big)= \sum_{d_A+d_B=d} (A\cup B, d_A, d_B).$$
\begin{lm}
\label{pb}
The relations among the boundary components $\deli'$ in $Pic(\M)\otimes \Q$
are the $\eta$ pull-backs of the relations among the boundary
components of $\barr{M}_{0,n}$.
\end{lm}
\begin{pf}
Let $\pi:S=\proj^1\times C \rarr C$ be a family with
$n$ sections. Let $\mu: S \ - \rarr \proj^r$ be a rational
maps with simple base points obtained from a line bundle
of degree type $(d,k)$. Suppose the special points satisfy (1), (2),
and
\begin{enumerate}
\item[($3'$.)] An intersection point lies on at most $n-2$ sections.
\item[(4.)] Every simple base point is an intersection point.
\end{enumerate}
Note condition ($3'$) implies condition (3).
For large $k$, the simple base points may be selected
arbitrarily (with arbitrary degree) among the intersection points.
Let $\barr{S}$ be the blow-up of $S$ at the special points; let
$\lambda: C\rarr \M$ be the induced curve.
By condition ($3'$), the
family $\barr{S} \rarr C$ with the strict transforms of the
sections is flat family of stable, $n$-pointed, genus $0$ curves.
The induced morphism $\gamma: C \rarr \barr{M}_{0,n}$ is simply
$\gamma= \eta\circ \lambda$.

Suppose $\sum_{K\in \deli'} c_K K =0$ is a relation
in $Pic(\M)\otimes \Q$ ($c_B \in \Q$).
Let $$K=(A\cup B, d_A,d_B)\in \deli'.$$ Let $(A\cup B)$ be corresponding
boundary component of $\barr{M}_{0,n}$.  The set
theoretic intersection $C\cdot K$ is the subset
of $C \cdot  (A\cup B)$ with simple base points of the correct degree.
Since the simple base points can be assigned arbitrary degrees,
the coefficient $c_K$ must depend only on the partition $(A\cup B)$
and not on the weights $d_A, d_B$. It now follows the relation
$\sum_{K\in \deli'} c_K \cdot K =0$ must be the $\eta$ pull-back
of a boundary relation in $\barr{M}_{0,n}$.
\end{pf}
In particular, it follows there are ${n-1 \choose 2}-1$ independent
relations among the boundary components $\deli'$. For $n\geq 4$,
$$|\deli|= d+ dn+ |\deli'|,$$
$$|\deli'|= (d+1) \cdot (2^{n-1}-1-n).$$
By Lemmas (\ref{n3}), (\ref{ihh}), (\ref{i01}), (\ref{pb}),
the Picard number of $\M$  ($n\geq4$) is:
$$dim\ Pic(\M)\otimes \Q = (d+1)\cdot 2^{n-1}-{n\choose 2}.$$
All the numerical relations are obtained from
linear equivalences.
The proof of Proposition (\ref{dimm}) is complete.

\section{Computations in $Pic(\M) \otimes \Q$}
\label{calc}
\subsection{The Universal Curve and $\pi_*(c_1(\omega_{\pi})^2)$}
\label{wclass}
Classes of certain canonical elements in $Pic(\M)\otimes \Q$ will be computed
via intersections with curves.
These computations will be used in the
proof of Proposition (\ref{top}).
In order to use the coarse moduli space throughout, an automorphism
result is required.
\begin{lm}
\label{auto}
Let $d\geq0$, $g=0$, $r\geq 2$. The locus of Kontsevich
stable maps in $\M_{0,n}(r,d)$ with nontrivial automorphisms is of codimension
at least $2$ except in one case: $\M_{0,0}(2,2)$.
\end{lm}
\begin{pf}
The assertion follows from naive dimension estimates.
If $d=0$ or $1$ , the are no stable maps with nontrivial automorphisms. Let
$\M=\M_{0,n}(r,d)$,  $(d\geq 2, r\geq 2)$. Recall
$dim \M = rd+d+r+n-3$. Certainly the generic elements
of the boundary components are automorphism-free.
Let $A\subset \M$ be the locus of non-boundary, stable maps
with nontrivial automorphisms.
If a map $\mu: \proj^1 \rarr \proj^r$ with $n$ distinct marked
points has an nontrivial automorphism, $\mu$ must be
a $k\geq 2$ to $1$ map. For fixed $2\leq k\leq d$,
the map $\mu$ moves in a family of dimension at most:
$$(r+1)\cdot ({d\over k}+1)-1-3 + 2\cdot (k+1)-1-3=
(rd+d)\cdot {1\over k}+r-3+2k-2.$$
The $n$ marked points must be fixed points of the nontrivial automorphism
and hence move in a zero dimensional family for each $\mu$.
A calculation yields:
\begin{eqnarray*}
dim \M- dim A & \geq & (rd+d)\cdot (1- {1\over k}) +n -2k+2 \\
& = & rd+d+n+2 - {rd+d+2k^2\over k}.
\end{eqnarray*}
A study of the function $(rd+d+2k^2)/k$ for $2\leq k \leq d$ shows the
maximum value must be
attained at the end points $k=2, d$.
If $k=2$,
$$ rd+d -{rd+d+8\over 2} = (r+1){d\over 2}-4 \geq 0$$
except when $r=2$, $d=2$.
If $k=d$,
$$ rd+d-{rd+d+2d^2\over d}= (r-1)(d-1)-2  \geq 0 $$
except when $r=2$, $d=2$.
$A$ is of codimension at least 2 all cases except $\M_{0,0}(2,2)$.
\end{pf}

Since $\M_{0,0}(2,2)$ is isomorphic to the space of complete
conics, its intersection theory is well known. In the sequel,
it will be assumed $(g,n,r,d)\neq (0,0,2,2)$.
Let $\M^*\subset \M$ denote the automorphism-free locus.
There is a universal Kontsevich stable family of maps over
$\M^*$:
$$\pi: U^* \rarr \M^*$$
with sections $s_1, s_2 \ldots, s_n$ and a morphism
$$\mu: U^* \rarr \proj^r.$$
See [P] for details. Let $\omega_{\pi}$ be the
relative dualizing sheaf of $\pi$.
Since the complement of $\M^*$ is of codimension at
least $2$ in $\M$, the following are well-defined
elements of $Pic(\M)\otimes \Q$ :
\begin{equation}
\label{fcla}
\pi_*( c_1(\omega_{\pi})^2), \ \ \pi_*( s_i^2).
\end{equation}
Since $Pic(\M) \otimes \Q \eqq Num(\M) \otimes \Q$, explicit expressions
of the classes (\ref{fcla}) in terms of the generators
$\{\LL_i\} \cup \deli \cup \{\HH\}$ can be found by calculating intersection
products
with curves in $\M$. The methods of section (\ref{realy}) will be used
to determine curves in $\M$.
First consider $\pi_*(c_1(\omega_{\pi})^2)$:

\begin{lm}
\label{om} For $d\geq 0$, $g=0$, $r\geq 2$,
$\ \pi_*( c_1(\omega_{\pi})^2) = - \sum_{K\in \deli} K$ in
$Pic(\M)\otimes \Q$.
\end{lm}
\begin{pf}
Let $\pi:S\rarr C$ be a projective bundle of rank $1$ over
a nonsingular curve $C$. Let $\omega_{\pi}$ be the
relative dualizing sheaf. A simple computation yields
$$\pi_*(c_1(\omega_{\pi})^2)= 0$$
in $Num(C)$.
Let $\rho: \barr{S}\rarr S$ be the blow-up at $k$ points
in distinct fibers of $\pi$. Let $\barr{\pi}:\barr{S} \rarr C$
be the composition.
$$\omega_{\barr{\pi}}= \rho^*(\omega_{\pi}) + \sum_{i=1}^{k} E_i$$
where the $E_i$ are the exceptional divisors of $\rho$. Hence,
$$\pi_*(c_1(\omega_{\barr{\pi}})^2)= -k$$
in $Num(C)$.
By considering curves $C\rarr \M^*$ and the pull-back of $U^*$,
it follows $\pi_*( c_1(\omega_{\pi})^2) = -\sum_{K\in \deli} K$
in $Num(\M) \otimes \Q$. By Proposition (\ref{dimm}), the
Lemma is proven.
\end{pf}

\subsection{The Class $\pi_*(s_1^2)$}
\label{sclass}
The determination of the class $\pi_*(s_1^2)$ is
surprisingly different in the cases $d=0$ and $d\geq 1$.
If $d=0$, it suffices to determine $\pi_*(s_1^2)$
for the universal family over $\barr{M}_{0,n}$.
Let $\deli$ be the set of boundary components of $\barr{M}_{0,n}$.
There is a partition of $\deli$ with respect to the first
marking. For $2\leq j \leq n-2$, let $\deli^{1}_j \subset \deli$
be defined by:
$$(A\cup B) \in \deli^{1}_j \ \ \ if \ and\ only \ if \ \ \ 1\in A, \ |A|=j.$$
There is a disjoint union
$$\deli= \bigcup_{j=2}^{n-2} \deli^{1}_j.$$
Let $K^{1}_{j}=\sum_{K\in \deli^{1}_{j}} K$.
\begin{lm}
\label{self0} The class
$\pi_*(s_1^2)$ is expressed in $Pic(\barr{M}_{0,n})\otimes \Q$ by:
\begin{equation}
\label{exx0}
\pi_*(s_1^2)=
-{1\over {n-1\choose 2}} \cdot \sum_{j=2}^{n-2} {n-j\choose 2} K^{1}_{j}.
\end{equation}
\end{lm}
\begin{pf}
The proof is by intersections with curves in $\barr{M}_{0,n}$.
Let $S=\proj^1\times C$ be a family with $n$ sections $s_1, \ldots, s_n$.
Let $s_1$ be of degree type $(1,q)$. For $2\leq i \leq n$,
let $s_i$ be of type $(1,p_i)$. As usual, assume the blow-up $\barr{S}$
up of $S$ at the intersection points yields a family of stable, $n$-pointed
curves over $C$ with at most one exceptional divisor in each fiber.
Let $\lambda: C \rarr \barr{M}_{0,n}$ be the induced map. It will be checked
that the left and right sides of (\ref{exx0}) have the same intersection with
$C$.

A point of $C\cdot K^{1}_{j}$ can arise in exactly two cases. First an
intersection point of $j$ sections including $s_1$ can be blown-up.
Second, an intersection point of $n-j$ sections not including $s_1$ can
be blown-up. Let
$$C\cdot K^{1}_{j}= x_j + y_j$$
where $x_j$, $y_j$ are the number of instances of the first and
second cases respectively. Let $\barr{s}_1$ be the strict transform
of $s_1$ in $\barr{S}$. The intersection of $C$ with the left side of
(\ref{exx0}) is :
$$\barr{\pi}_*(\barr{s}_1^2)= 2q- \sum_{j=2}^{n-2} x_j.$$
For $2\leq  i \leq n$, $s_i$ intersects $s_1$ in $q+p_i$ points.
The following equation is easily obtained by analyzing intersection
points contained in $s_1$:
\begin{equation}
\label{ff}
(n-1) q + \sum_{i=2}^{n} p_i = \sum_{j=2}^{n-2} (j-1) x_j.
\end{equation}
Similarly, the number of intersections of the sections $2\leq i \leq n$
among themselves is $(n-2)\cdot \sum_{i=2}^{n} p_i$. Analysis of
intersection points not contained in $s_1$ yields:
\begin{equation}
\label{ss}
(n-2) \cdot \sum_{i=2}^{n} p_i = \sum_{j=2}^{n-2}
{j-1 \choose 2} x_j + {n-j\choose 2} y_j.
\end{equation}
Via equations (\ref{ff}) and (\ref{ss}),
\begin{eqnarray*}
{n-1\choose 2} \cdot (2q-\sum_{j=2}^{n-2} x_j)& = &
\sum_{j=2}^{n-2} \big((n-2)(j-1)-{j-1\choose 2}-{n-1\choose 2}\big) x_j
- {n-j\choose 2} y_j
\\
& = & - \sum_{j=2}^{n-2} {n-j\choose 2} (x_j+y_j).
\end{eqnarray*}
The Lemma is proved.
\end{pf}

Consider now the case $d\geq 1$.
Let $\M=\M_{0,n}(r,d)$ where $d\geq1$, $n\geq 1$. Let $1$ be the first
marking. There is another partition of $\deli$ with respect to the
first marking depending upon the degree. For $0\leq j \leq d$, let
$\deli^{1,j}\subset \deli$ be defined by:
$$(A\cup B, d_A, d_B) \in \deli^{1,j} \ \ \
if\ and\ only\ if \ \ \ 1\in A,\ d_A=j.$$
Note if $n=1$, then $\deli^{1,0}, \deli^{1,d}=\emptyset$.
If $n=2$, $\deli^{1,d}=\emptyset$.
In all other cases $\deli^{1,j}\neq \emptyset$.
There is a disjoint union
$$\deli= \bigcup_{j=0}^{d} \deli^{1,j}.$$
Let $K^{1,j}= \sum_{K\in \deli^{1,j}} K$. Let $K^{1,j}=0$ if
$\deli^{1,j}=\emptyset$.
\begin{lm}
\label{self} In case $d\geq 1$,
The class
$\pi_*(s_1^2)$ is expressed in $Pic(\M)\otimes \Q$ by:
\begin{equation}
\label{exx}
\pi_*(s_1^2)= -{1\over d^2}\HH + {2\over d} \LL_1
-\sum_{j=0}^{d} {(d-j)^2\over d^2} K^{1,j}.
\end{equation}
\end{lm}
\begin{pf}
The proof is by intersections with curves in $\M$.
Let $\pi: S=\proj^1\times C \rarr \proj^1$ be a family
with $n$ sections $s_1, \ldots, s_n$.
Let $s_1$ be of degree type $(1,q)$.
Let $\mu:S - \ \rarr \proj^r$
be a rational map with simple base points obtained
from a line bundle of degree type $(d,k)$. Let conditions
(1), (2), ($3'$), (4) of section (\ref{realy})
be satisfied. Let $\barr{S}\rarr S$ be the
blow-up at the special points. Let $\lambda: C \rarr \M$
be the induced map. It will be checked that the left and right sides
of (\ref{exx}) have the same intersection with $C$.

A point of $C\cdot K^{1,j}$ can arise in exactly two cases. First,
a simple base point of degree $j$ contained in $s_1$ can be blown-up.
Second, a simple base point of degree $d-j$ not contained in $s_1$ can
be blown-up. Let
$$C \cdot K^{1,j}= x_j + y_j$$
where
$x_j$, $y_j$ are the number of instances of the first and
second cases respectively. Let $\barr{s}_1$ be the strict
transform of the section $s_1$ to $\barr{S}$. The intersection
of $C$ with the left side of (\ref{exx}) is given by:
$$\barr{\pi}_*(\barr{s}_1^2) = 2q- \sum_{j=0}^{d} x_j.$$
A straightforward computation yields:
$$C \cdot \HH= 2dk - \sum_{j=0}^{d} j^2 x_j - \sum_{j=0}^{d} (d-j)^2 y_j,$$
$$C \cdot \LL_1= dq+k - \sum_{j=0}^{d} j x_j.$$
The equality of the intersection of $C$ with the left and
right sides of (\ref{exx}) is now a matter of simple algebra.
\end{pf}

\subsection{The Class $\TT$}
\label{tclass}
Let $\M=\M_{0,n}(r,d)$, $d\geq 2$. Let $H\subset \proj^r$ be a hyperplane.
A tangency Weil divisor $\TT_H \subset \M$ is defined as follows.
Let $W_H\subset \M$ be the open locus of maps
$\mu:C \rarr \proj^r$ where $\mu^{-1}(H)$ is a
subscheme of $d$ reduced points of $C_{nonsing}$.
Let $\TT_H$ be the complement of $W_H$.

It must be shown that $\TT_H$ is of pure codimension $1$ in $\M$.
Let $M_H\subset \M$ be the open locus of maps $\mu: C \rarr \proj^r$
satisfying :
$$\forall x\in\mu^{-1}(H), \ \ \
x\in C_{nonsing} \ \ and \ \ d\mu_x\neq 0.$$ The intersection
$\TT_H\cap M_H$ corresponds to geometric tangencies and is certainly
of pure codimension $1$ in $M_H$ ($d\geq 2$). The complement
$\M \setminus M_H$ is of codimension $2$ in $\M$.
It is not hard to
see the closure of $\TT_H\cap M_H$ in $\M$
contains the complement $\M \setminus M_H$. Therefore, $\TT_H$ is
a Weil divisor.

Define
for $0\leq j \leq [{d\over2}]$, $\deli^j \subset \deli$ as follows. A boundary
component $(A\cup B, d_A, d_B)\in \deli^j$ if and only if the
degree partition $d_A+d_B=d$ equals the partition $j+(d-j)=d$.
Let $K^j= \sum_{K\in \deli^j} K$.
\begin{lm}
\label{tan}
The class of $\TT$ can be expressed in $Pic(\M)\otimes \Q$ by:
\begin{equation}
\label{vexx}
\TT = {d-1\over d} \HH + \sum_{j=0}^{[{d\over2}]} {j(d-j)\over d} K^j.
\end{equation}
\end{lm}
\begin{pf}
Let $S$, $\mu$, $\barr{S}$, $\lambda: C \rarr \M$ be
exactly as in the proof of Lemma (\ref{self}).
It will be checked that the left and right sides of
(\ref{vexx}) have the same intersection with $C$.

As before, a point of the intersection $C\cdot K^j$ can arise
in two cases. A simple point of degree $j$ or $d-j$ can be blown-up.
Let $$C\cdot K^j = x_j + y_j$$
where $x_j$ and  $y_j$ are the number instances of the first and
second case respectively. Let $E_{x_j}$ be the union of the
$x_j$ exceptional divisors in $\barr{S}$ obtained from the
$x_j$ points of $C\cdot K^j$. Let $E_{y_j}$ be defined
similarly.

First, the intersection $C\cdot \TT$ is calculated.
A general element of $\barr{\mu}^*(\oh_{\proj}(1))$ is
a nonsingular curve $D$ in the linear series $(d,k)-
\sum_{j} j E_{x_j} - \sum_{j} (d-j) E_{y_j}$. Adjunction yields:
$$2 g_D-2 = d(2g_C-2)+ 2dk-2k- \sum_{j=0}^{[{d\over2}]} j(j-1) x_j -
\sum_{j=0}^{[{d\over2}]} (d-j)(d-j-1) y_j.$$
Since $D$ is a $d$ sheeted cover of $C$, the Riemann-Hurwitz formula
determines the ramifications:
$$C\cdot \TT = 2dk-2k - \sum_{j=0}^{[{d\over 2}]} j(j-1)x_j+(d-j)(d-j-1)y_j.$$
$C\cdot \HH$ is simply $D^2$. Hence
$$C\cdot \HH= 2dk - \sum_{j=0}^{[{d\over 2}]} j^2x_j+ (d-j)^2 y_j.$$
Again, an algebraic computation yields the equality of the
intersections of the left and right sides with $C$.
\end{pf}

\section{Intersections of $\Bbb{Q}$-divisors}
\label{inter}
\subsection{Intersections of the classes $\{\LL_i\}\cup \{\HH\}$}
The top dimensional intersection products on $\M=\M_{0,n}(r,d)$ of the
classes $\{\LL_i\}$ are algorithmically determined by the
First Reconstruction Theorem [K-M]. These top
classes are computed recursively in $d$ and $n$.
The algorithm requires one initial value: the number
of lines in $\proj^r$ through two points. The top intersection
products of $\{\LL_i\}$ are exactly the characteristic numbers ($\beta=0$) of
rational curves in $\proj^r$.

Top dimensional intersections of the classes $\{\LL_i\} \cup \{\HH\}$
are also characteristic numbers of rational curves in $\proj^r$.
Each factor of $\HH$ is a codimension-2  characteristic condition. For example,
if $\M=\M_{0,0}(2,3)$, then $\HH^8$ equals the number of
rational plane cubics through $8$ general points. If $\M=\M_{0,2}(3,4)$,
then $c(\LL_1)^3\cdot c(\LL_2)^3 \cdot \HH^{12}$ equals the number
of rational space quartics passing through  2 general points and meeting
12 general lines.

\subsection{Boundary Components}
Let $K=(A\cup B, d_A, d_B)$ be a boundary component of $\M_{0,n}(r,d)$.
Let $\M_A=\M_{0, |A|+1}(r,d_A)$ and $\M_B= \M_{0, |B|+1}(r, d_B)$.
Let the additional markings be $p_A$ and $p_B$ respectively.
Let $e_A: \M_A \rarr \proj^r$ and $e_B: \M_B\rarr \proj^r$ be the
evaluation maps obtained from the markings $p_A$ and $p_B$.
Let $\tau_A$, $\tau_B$ be the projections from $\M_A \times \M_B$
to the first and second factors respectively.
Let $\tl {K}= \M_A \times_{\proj^r} \M_B$ be the fiber
product with respect to the evaluation maps $e_A$, $e_B$.
$\tl{K}\subset \M_A \times_{\com} \M_B$ is the closed
subvariety
$(e_A\times_{\com} e_B)^{-1} (D)$ where $D\subset \proj^r \times \proj^r$ is
the diagonal. $\tl{K}$ is easily seen to be an irreducible, normal,
projective variety with finite quotient singularities. These results
follow, for example, from the local construction given in [P].
The class of $\tl{K}$ in $\M_A \times \M_B$ can be  computed by the
pull-back of the K\"unneth
decomposition of the diagonal in $\proj^r \times \proj^r$:
\begin{equation}
\label{f1}
[\tl{K}] = \sum_{i=0}^{r} \tau_A^*(c_1(\LL_A)^i)\cdot
\tau_B^*(c_1(\LL_B)^{r-i})
\end{equation}
where $\LL_A$, $\LL_B$ are the line bundles on $\M_A$, $\M_B$ induced
by the marking $p_A$, $p_B$.

There is a natural map $\psi: \tl{K} \rarr K$.
The set theoretic description of $\psi$ is clear:
$\psi([\mu_A],[\mu_B])$
is the moduli point of the map obtained by gluing maps $\mu_A$, $\mu_B$
along the markings
$p_A$, $p_B$. It is not hard to define $\psi$ algebraically.
$\psi$ is a birational morphism
except when $n=0$ and $d_A=d_B=d/2$. In the latter case, $\psi$ is generically
2-1.

The pull-backs of the classes
$\{\LL_i\} \cup \deli \cup \{\HH \}$ on $\M$ to $\tl{K}$ are determined in
the following manner. Let $\HH_A$, $\HH_B$ be the codimension-$2$ plane
incidence classes on $\M_A$, $\M_B$. Clearly,
\begin{equation}
\label{f2}
\psi^*(\HH)= (\tau_A^*(\HH_A) + \tau_B^*(\HH_B)) \ |_{\tl{K}}.
\end{equation}
Let $P$ be the marking set of $\M$. For each $i\in P$, $i$ is either
in $A$ or $B$. It follows
\begin{equation}
\label{f3}
\psi^*(\LL_i)= \tau_A^*(\LL_i) \ |_{\tl{K}}, \ \ \psi^*(\LL_i)=
\tau_B^*(\LL_i)\ |_{\tl{K}}
\end{equation}
in case $i\in A$, $i\in B$ respectively.

Let $T=(A'\cup B', d_{A'}, d_{B'})$ be a boundary component of $\M$ {\em not}
equal to $K$. $T$  intersects
$K$ exactly when one of the following two conditions hold:
\begin{enumerate}
\item[(i.)] There exists a subset $C\subset A$ and an integer $d_C$
such that $$( (A\setminus C) \cup (B\cup C), d_A-d_C, d_B+d_C)= T.$$
\item[(ii.)] There exists a subset $C\subset B$ and an integer $d_C$
such that $$( (A\cup C) \cup (B\setminus C), d_A+d_C, d_B-d_C)=T.$$
\end{enumerate}
\begin{equation}
\label{f4}
\psi^*(T) =
\sum_{C, d_c} \tau_A^*(A\cup (C\cup \{p_A\}), d_A, d_C) \ |_{\tl{K}}
\end{equation}
$$ + \ \ \sum_{C,d_c} \tau_B^* (B\cup (C\cup \{p_B\}), d_B, d_C) \
|_{\tl{K}}.$$
The sums on the right are taken over subsets $C$ and degrees $d_c$
that satisfy (i) and (ii) above respectively. The main point is
distinct boundary divisors have transverse (if nonempty) intersections
in the stack $\barr{\cal{M}}_{0,n}(r,d)$. This can be seen an a property
inherited from the Mumford-Knudsen space $\M_{0,m}$ by the local
construction given in [P]. Since the automorphism loci of
$\M_{0,n}(r,d)$ and the boundary component $(A\cup B, d_A, d_B)$
are of codimension at
least two in $\M_{0,n}$, $(A\cup B, d_A, d_B)$ respectively,
the transverse intersection property descends to
the coarse moduli space.

Let $\omega_{\pi A}$, $\omega_{\pi B}$ denote the relative dualizing
sheaves of the the universal families over $\M^*_A$, $\M^*_B$ respectively.
There are two universal curves over $\tl{K}^*=\tl{K}\cap (\M^*_A\times \M^*_B)$
obtained via pull-back of the universal families
$U_A^*$ and $U_B^*$. These universal curves glue on the
sections $s_{pA}$ and $s_{pB}$ to form a universal family
$$\tl{\pi}: U^*_{\tl{K}^*} \rarr \tl{K}^*$$ of maps for the
moduli problem of $\M$.
It follows,
$$\omega_{U^*_{\tl{K}^*}}\ |_ {\tau_A^*(U_A^*)}=
\tau_A^*(\omega_{\pi A}) + s_{pA},$$
$$\omega_{U^*_{\tl{K}^*}}\ |_ {\tau_B^*(U_B^*)}=
\tau_B^*(\omega_{\pi B}) + s_{pB}.$$
Hence
$$\psi^*(\pi_*(c_1(\omega_{\pi})^2))= \tau_A^*(
\pi_{A*}((c_1(\omega_{\pi A})+s_{pA})^2))
+ \tau_B^*(\pi_{B*}((c_1(\omega_{\pi B})+s_{pB})^2)).$$
A normal bundle calculation yields $c_1(\omega_{\pi A})\cdot s_{pA}=
-s_{pA}^2$.
Hence,
$$(c_1(\omega_{\pi A})+s_{pA})^2 = c_1(\omega_{\pi A})^2 - s_{pA}^2$$
(similarly for $B$).
Recall $\pi_*(c_1(\omega_{\pi})^2)= -\sum_{T\in \deli} T$.
Finally,
\begin{equation}
\label{f5}
-\psi^*(K)  = \sum_{T\in \deli, T\neq K} \psi^*(T) \ \
+\ \ \tau_A^*(\pi_{A*}(c_1(\omega_{\pi A})^2 - s_{pA}^2))
\end{equation}
$$+ \ \ \tau_B^*(\pi_{B*}(c_1(\omega_{\pi B})^2 - s_{pB}^2)).$$
Lemmas (\ref{om}), (\ref{self0}), and (\ref{self}) express
$\pi_{A*}(c_1(\omega_{\pi A})^2)$, $\pi_{A*}(s_{pA}^2)$ explicitly
in terms of the standard classes $\{\LL_i\} \cup \deli \cup \{\HH\}$
on $\M_A$ (similarly for $\M_B$). Via equations (\ref{f2}) - (\ref{f5}),
the $\psi$ pull-back of every standard class $\{\LL_i\} \cup \deli \cup
\{\HH\}$ on $\M$ has now been expressed as the restriction to
$\tl{K}$ of a  linear combination of the
$\tau_A$ and $\tau_B$ pull-backs of standard classes on $\M_A$ and
$\M_B$.

\subsection{The Algorithm}
The inductive algorithm for computing top intersection products is
now clear. All top monomials in the elements $\{\LL_i\}\cup \{\HH\}$
are known by the
First Reconstruction theorem. If a monomial product on $\M$ includes
a boundary class $K$, the intersection is carried out on $\tl{K}$.
By the above formulas (\ref{f1})-(\ref{f5}), the desired monomial can be
expressed as a sum of top products of standard classes
on $\M_A$ and $\M_B$. Since $\M_A$ is of lesser degree or of lesser marking
number than $\M$ (similarly for $\M_B$), the inductive process terminates.

\subsection{Characteristic Numbers}
\label{cnum}
Lemma (\ref{tan}) expresses the hyperplane tangency condition in terms
of the standard classes. Hence all top products of
the classes $\{\LL_i\} \cup \{\HH, \TT\}$ can be effectively
computed by the above algorithm.

It remains to check the top intersections of $\{\LL_i\} \cup \{\HH, \TT\}$
are the characteristic numbers of rational curves. Let
\begin{equation}
\label{cy}
c_1(\LL_1)^{l_1} \cdots c_1(\LL_n)^{l_n} \cdot \HH^{\alpha}\cdot \TT^{\beta}
\end{equation}
be a top product on $\barr{M}=\barr{M}_{0,n}(r,d)$.
Since the $\LL_i$ are pull-backs of $\oh_{\proj^r}(1)$ via the evaluation
maps, codimension $l_i$ linear spaces of $\proj^r$ determine
representatives of $c_1(\LL_i)^{l_i}$. The cycle $\HH^{\alpha}$ is
determined by $\alpha$ codimension $2$ linear spaces in $\proj^r$. Finally,
the cycle $\TT^{\beta}$ is determined by $\beta$ hyperplanes in $\proj^r$.
When ($\beta\neq 0$), it is assumed $d\geq 2$.
The first step
is to show for general choices of all the linear spaces of $\proj^r$
in question,
the intersection cycle (\ref{cy}) in $\barr{M}$ is at most 0  dimensional and
corresponds (set theoretically) to the correct geometric locus. The second step
is to show the intersection cycle is multiplicity free.

Let $\proj^{r*}$ be the parameter space of hyperplanes in $\proj^r$.
Defined the universal tangency subvariety
$$\TT_{univ} \subset \barr{M} \times \proj^{r*}$$
as follows. Let $W_{univ}\subset \M \times \proj^{r*}$ be the open locus of
pairs
$(\mu:C\rarr \proj^r, H)$ where $\mu^{-1}(H)$ is a subscheme
of $d$ reduced points of $C_{nonsing}$.
Let $\TT_{univ}$ be the complement of $W_{univ}$.
Let $\TT_H$ be the
the fiber of $\TT_{univ}$ over the parameter point of the
hyperplane $H$. $\TT_H$ is exactly the tangency Weil divisor
defined in section (\ref{calc}). Similarly, let
$$\HH_{univ} \subset \M \times  \grass(\proj^{r-2}, \proj^r)$$
be the universal codimension 2 plane incidence subvariety. The
fiber of $\HH_{univ}$ over the parameter point of the
codimension 2 plane $P$ is $\HH_P$.
Let $$I_{univ} \subset \M \ \times\  \grass(r-l_1,r) \times \cdots \times
\grass(r-l_n,r)\ \times\
\grass(r-2,r) \times \cdots \times  \grass(r-2,r)$$
$$\times \  \proj^{r*}\times \cdots \times \proj^{r*}$$
be the universal intersection cycle (\ref{cy})
defined by the universal divisors $\TT_{univ}$, $\HH_{univ}$ and
the evaluation maps. $I_{univ}$ is closed subvariety.

In the first step, slightly more than the dimensionality of
the general intersection cycle will be established.
A map $\mu: C \rarr \proj^{r}$ is {\em simply tangent} to
a hyperplane $H$ if
\begin{enumerate}
\item[(i.)] $\mu^{-1}(H)\subset C_{nonsing}$.
\item[(ii.)] As a subscheme, $\mu^{-1}(H)$ consists of $1$ double and $d-2$
reduced points.
\end{enumerate}
A map $\mu: C \rarr \proj^{r}$ has {\em simple intersection} with
a codimension 2 plane $P$ if
\begin{enumerate}
\item[(i.)] $\mu^{-1}(P)$ consists of 1 point $x \in C_{nonsing}$.
\item[(ii.)] $Im(d\mu(x))$ and the tangent space of $P$ span
maximal rank.
\end{enumerate}

\begin{lm}
For general choices of linear spaces
\begin{equation}
\label{lin}
L_1, \ldots, L_n,\  P_1, \ldots, P_{\alpha},\  H_1, \ldots, H_{\beta}
\end{equation}
the intersection
cycle (\ref{cy}) is at most 0 dimensional and set theoretically
corresponds to maps $\mu: C \rarr \proj^r$ satisfying:
\begin{enumerate}
\item[(1.)] $C\eqq \proj^1$, $\mu$ is an immersion/embedding ($r=2$ /
$r\geq 3$).
\item[(2.)] $\forall k$, $\mu$ is simply tangent to the hyperplanes $H_k$.
\item[(3.)] $\forall j$, $\mu$ intersects the linear spaces $P_j$ simply.
\item[(4.)] $\forall i$, the $\mu$-image of the $i^{th}$ marked point
lies in $L_i$.
\end{enumerate}
\end{lm}
\begin{pf}
The intersection cycle $I$ determined by the
linear spaces (\ref{lin}) is the fiber of $I_{univ}$ over
the parameter points of the linear spaces.
$dim(I)\leq 0$ is an open condition in the parameter
space. It is first checked that general choice of the linear spaces (\ref{lin})
yields an intersection cycle of dimension at most $0$.

Let $[\mu]\in \M$ be the moduli point of a map $\mu: C \rarr \proj^r$.
By Bertini's Theorem, the general hyperplane $H$ is transverse to $\mu$.
Therefore, the general tangency divisor $\TT_H$ satisfies
$[\mu] \notin \TT_H$. Similarly, the general incidence divisor
$\HH_P$ satisfies $[\mu]\notin \HH_P$.
By choosing at each stage tangency and incidence divisors
that reduce the dimension of every component of the intersection,
it follows
$$\HH_{P_1} \cap \ldots \cap \HH_{P_\alpha} \cap
\TT_{H_1} \cap \ldots \cap \TT_{H_\beta}$$
has codimension at least $\alpha + \beta$.
Since the remaining intersections are obtained from basepoint
free linear series, the general intersection cycle has
dimension at most $0$.

If the general parameter point yields an empty cycle $I$,
there is nothing more to prove. Let $W$ be the open set
of the parameter space where $dim(I)=0$.
The conditions (1-3) on $I$ determine open
sets $W_1, W_2, W_3 \subset W$. Condition (iv) is automatic. It suffices
to show $W_i$ is nonempty for $1\leq i \leq 3$.

The subset $Y\subset \M$ of maps that are not immersion/embedding
($r=2$ / $r\geq 3$) is of codimension at least $1$. Hence, by the dimension
reduction argument above, $Y\cap I=\emptyset$ for a general parameter
point. Therefore, $W_1 \neq \emptyset$.

Let $W_{2,k}, W_{3,j} \subset W$,
be the set of parameter points that satisfy
condition (2), (3)  for the hyperplane $H_k$,
linear space $P_j$ respectively. Since $W_2= \cap _{k=1}^{\beta} W_{2,k}$ and
$W_3=\cap_{j=1}^{\alpha} W_{3,j}$,
it suffices to show $W_{2,k}, W_{3,j} \neq \emptyset$. Let
$H_k$ be any hyperplane. The locus of moduli points
$[\mu]\in \TT_{H_k}$ that are not simply tangent is of codimension at least
2 in $\M$. By the dimension reduction argument, $W_{2,k}\neq \emptyset$.
Similarly, the locus of moduli points $[\mu]\in \HH_{P_j}$
that do not intersect simply is of codimension at least 2 in $\M$.
As before $W_{3,j} \neq \emptyset$.
\end{pf}

It must now be shown that the intersection cycle (\ref{cy}) is reduced for
general linear spaces. This transversality is established by
Kleiman's Bertini Theorem. Unfortunately, since
the divisors $\TT_H$, $\HH_P$ need not move {\em linearly}, Bertini's Theorem
can not be directly applied to $\barr{M}$. Instead, an auxiliary
construction is undertaken. Kleiman's Bertini Theorem is applied to
the universal curve over $\barr{M}$. It will be shown that suitable
transversality on the universal curve implies transversality on $\barr{M}$.

Let $\barr{M}^0\subset \barr{M}$ be
the open set of immersed/embedded ($r=2$, $r\geq 3$) maps with
irreducible domains. Since (for general linear spaces) the
intersection cycle (\ref{cy}) lies in $\barr{M}^0$, transversality need only be
established in $\barr{M}^0$. Note $\barr{M}^0$ is in the automorphism-free
locus. Let $U \rarr \barr{M}^0$ be the universal curve.
Let $\mu: U \rarr \proj^r$ be the universal map. $U$, $\barr{M}^0$ are
nonsingular.
Let $\proj T$ be the projective tangent bundle of $\proj^r$.
Since each point of $\barr{M}^0$ corresponds to an immersion/embedding,
there is a natural algebraic map
$\nu: U \rarr \proj T$ given by the differential of $\mu$. The
map $\nu$ is a  lifting of $\mu$.

By projectivizing tangent spaces,
the hyperplanes $H_1, \ldots, H_{\beta}$ define nonsingular, codimension
2 subvarieties of $\proj T$: $$\proj H_1, \ldots, \proj H_{\beta}$$
Let $U_1, \ldots, U_{\beta}$ be $\beta$
copies of the universal curve $U$. Let $U'_1, \ldots, U'_{\alpha}$
be $\alpha$ more copies of $U$. Define the product:
$$X \eqq U'_1 \times_{\barr{M}^0} \ldots
\times_{\barr{M}^0} U'_{\alpha} \times_{\barr{M}^0} U_1 \times_{\barr{M}^0}
\ldots \times_{\barr{M}^0} U_{\beta}.$$
Let $\mu'_j: X \rarr \proj^r$, $\nu_k:X \rarr \proj T$ be the
maps obtained by projection onto $U'_j$, $U_k$ and composition with $\mu$,
$\nu$
respectively.

Kleiman's Bertini Theorem may now be applied. The group $GL_{r+1}(\com)$ acts
transitively on $\proj^r$, $\proj T$. Hence, the general intersection
$$\mu'_1\ ^{-1}(P_1)\cap \ldots \cap \mu'_{\alpha}\ ^{-1}(P_{\alpha})
\cap \nu_1^{-1}(\proj H_1) \cap \ldots \cap \nu_{\beta}^{-1}(\proj H_{\beta})
\subset X$$
is nonsingular and of the correct codimension (if nonempty).

It remains to obtain the corresponding result on $\barr{M}$.
Consider the nonsingular, codimension 2 subvariety
$\mu^{-1}(P_j) \subset U$. The projection $\mu^{-1}(P_j) \rarr \HH_{P_j}
\cap \barr{M}^0$ is \'etale and 1-1 over the locus of of
maps meeting $P_j$ simply. Similarly, the projection $\nu^{-1}(\proj H_k)
\rarr \TT_{H_k}\cap \barr{M}^0$ is \'etale and 1-1 over the
the locus of maps simply tangent to $H_k$. From Lemma (\ref{sim}) below,
the projection
 $$\mu'_1\ ^{-1}(P_1)\cap \ldots \cap \mu'_{\alpha}\ ^{-1}(P_{\alpha})
\cap \nu_1^{-1}(\proj H_1) \cap \ldots \cap \nu_{\beta}^{-1}(\proj H_{\beta})
\longrightarrow\ \ \ \ $$
$$\ \ \ \ \ \ \HH_{P_1}\cap \ldots \HH_{P_\alpha} \cap \TT_{H_1} \cap \ldots
\cap \TT_{H_\beta} \cap \barr{M}^0$$
is \'etale and 1-1 over the locus of points in
$\HH_{P_1}\cap \ldots \HH_{P_\alpha} \cap \TT_{H_1} \cap \ldots
\cap \TT_{H_\beta} \cap \barr{M}^0$ corresponding to
 simple intersection and tangency.
It has therefore been proved, for general linear spaces, the
locus of $\HH_{P_1}\cap \ldots \HH_{P_\alpha} \cap \TT_{H_1} \cap \ldots
\cap \TT_{H_\beta} \cap \barr{M}^0$ corresponding to
 simple intersection and tangency is nonsingular and of the
correct codimension (if nonempty).
It was shown above, for general linear spaces, the intersection
cycle (\ref{cy}) involves only maps that have simple intersection
and tangency with the $P_j$, $H_k$.
Since the intersections $c_1(\LL_i)^{l_i}$ are obtained from
basepoint free linear series on $\barr{M}$, the further
intersections yield a reduced intersection cycle (\ref{cy}) by
Bertini's Theorem.

\begin{lm}
\label{sim}
Let $M$ be a nonsingular base. Let $\pi:U\rarr M$ be smooth
map of relative dimension 1. Let $D_1, D_2, \ldots, D_l\subset U$ be
nonsingular,
codimension $2$ subvarieties such that $D_i$ is \'etale and 1-1 over
$\pi(D_i)$.
Let $X\eqq U_1\times_M \ldots \times_M U_l$ be the fiber product of
copies of $U$. Let $\rho_i:X\rarr U_i$ be the projection. Then
$$\rho_1^{-1}(D_1) \cap \ldots \cap \rho_l^{-1}(D_l)\subset X$$
is \'etale and 1-1 over the intersection $\pi(D_1) \cap \ldots \cap \pi(D_l)
\subset M$.
\end{lm}
\begin{pf}
The issue is local on $M$.
Let $m\in M$ be in the intersection of the $\pi(D_i)$.
Choose local defining equations $(f_i)$ of
$\pi(D_i)$ near $m$. Let $u_i \in D_i$ be points over $m$.
Locally (in the analytic topology) at $u_i$, $U_i$ is an open set of the
trivial product
$\com_i \times M$
and $D_i$ is the intersection of $(f_i)$ with a section $(z_i)$ of this
product ($z_i$ is the coordinate on $\com_i$).
It now follows local equations for
for $\rho^{-1}_1 D_1\cap\ldots \cap \rho^{-1}_l D_l$ at $(u_1, \ldots, u_l)$
are
$(z_1, \ldots, z_l, f_1, \ldots, f_l)$ in
$\com_1 \times \ldots \com_l \times M$ which is
certainly \'etale over $(f_1,\ldots, f_l) \subset M$.
\end{pf}

All the characteristic numbers of rational
curves in projective space can be algorithmically computed. For example,
the number of twisted cubics in $\proj^3$ through $2$ points,
$6$ lines, and tangent to $2$ planes can be expressed as
$$c_1(\LL_1)^3 \cdot c_1(\LL_2)^3 \cdot
c_1(\LL_3)^2 \cdots c_1(\LL_8)^2 \cdot \TT^2$$
on $\M_{0,8}(3,3)$ or
$$c_1(\LL_1)^3\cdot c_1(\LL_2)^3 \cdot \HH^6\cdot \TT^2$$
on $\M_{0,2}(3,3)$.

\section{Examples}
\label{exam}
\subsection{Conics in $\proj^2$ and $\proj^3$}
Since the Hilbert schemes of lines and conics are Grassmanians and
projective bundles over Grassmanians, the $\beta=0$
characteristic numbers of
rational curves in degrees $1$ and $2$ can be calculated directly via
intersection theory on these Hilbert schemes. The tangency characteristic
numbers for conics classically required the beautiful space of complete
conics. $\barr{M}_{0,0}(2,2)$ is the space of complete conics. A
new calculation of the characteristic numbers for plane
conics is obtained by considering the pointed space
$\barr{M}_{0,1}(2,2)$.

Let $\M=\M_{0,1}(2,2)$. $Pic(\M)\otimes \Q$ is freely generated
by $\HH$, $\LL_1$, and the unique boundary component $K$
corresponding to the partition $(\{1\}\cup{\emptyset}, 1+1=2)$.
The top intersection numbers are ($dim\M_{0,1}(2,2)=6$):
$$\begin{array}{llrllllrllllr}
\HH^6 & & 0 & & & \HH^5K & & 0 & & & \HH^4K^2 & & 0 \\
\HH^5 \LL_1 & & +2 & & & \HH^4K  \LL_1 & & +6 & & &
\HH^3 K^2 \LL_1 & & +18\\
\HH^4 \LL_1^2 & & +1 & & & \HH^3K \LL_1^2 & & +3 & & &
\HH^2 K^2  \LL_1^2 & & +9 \\
\end{array}$$
$$\begin{array}{llrllllrllllrllllr}
\HH^3K^3&&0&&&\HH^2K^4&&0&&& \HH K^5 &&0&&&
K^6 &&0 \\
\HH^2 K^3\LL_1 &&-10 &&&
\HH K^4  \LL_1 && -30 &&& K^5 \LL_1 && +102 \\
\HH K^3 \LL_1^2 && -5 &&& K^4 \LL_1^2 && -15
\end{array}$$
\noindent Note $\LL_1^3=0$. The line tangency class $\TT={1\over2}(\HH+K)$
is determined by Lemma (\ref{tan}). The characteristc number
of plane conics through $\alpha$ points and tangent to
$\beta$ lines is ${1 \over 2}\HH^{\alpha}\TT^{\beta}\LL_1$:
$$\begin{array}{llr}
(1/ 2)\cdot\HH^5 \LL_1 & & 1 \\
(1/2)\cdot\HH^4 \TT \LL_1 & & 2 \\
(1/2)\cdot\HH^3 \TT^2 \LL_1 & & 4 \\
(1/2)\cdot\HH^2 \TT^3 \LL_1 & & 4 \\
(1/2)\cdot\HH \TT^4 \LL_1 & & 2 \\
(1/2)\cdot\TT^5  \LL_1 && 1
\end{array}$$
\noindent The class of
maps tangent to a fixed conic can be easily calculated by
the methods of Lemma (\ref{tan}). Let ${\cal{C}}\in Pic(\M)\otimes \Q$ denote
this conic tangency class. ${\cal{C}}=3\HH+K$.
The number of plane conics tangent
to 5 fixed conics is therefore ${1\over 2}{\cal{C}}^5\LL_1=3264$.

For
$r\geq 3$, $\barr{M}_{0,0}(r,2)$ differs from the space of complete
conics and the algorithm described above yields a new computation of
the characteristic numbers in these
cases.
Let $\M=\M_{0,0}(3,2)$. $Pic(\M)\otimes \Q$ is freely generated by
$\HH$ and the unique boundary component $K$ corresponding to the
degree partition $1+1=2$. $\tl{K}\subset \M_{0,1}(3,1) \times
\M_{0,1}(3,1)$. Since $\M_{0,1}(3,1)$ has no boundary, all top intersections
are known. Using the formulas of section (\ref{inter}), the
answers for the top intersections of $\HH$ and $K$ on
$\M_{0,0}(3,2)$  ($dim \M_{0,0}(3,2)=8$) are:
$$\begin{array}{llr}
\HH^8 & & +92 \\
\HH^7 K & & +140 \\
\HH^6  K^2 & & +140 \\
\HH^5  K^3 & & -100 \\
\HH^4  K^4 & & -68 \\
\HH^3  K^5 & & +172 \\
\HH^2  K^6 & & -20 \\
\HH K^7 & & -580 \\
K^8 & & +1820
\end{array}$$
By Lemma (\ref{tan}), $\TT= {1\over 2} (\HH+K)$. The
characteristic number of space conics through $\alpha$ lines
and tangent to $\beta$ planes is $\HH^\alpha  \TT^\beta$:
$$\begin{array}{llr}
\HH^8 & & 92 \\
\HH^7 \TT & & 116 \\
\HH^6  \TT^2 & & 128 \\
\HH^5  \TT^3 & & 104 \\
\HH^4 \TT^4 & & 64 \\
\HH^3  \TT^5 & & 32 \\
\HH^2 \TT^6 & & 16 \\
\HH  \TT^7 & &  8 \\
\TT^8 & & 4
\end{array}$$
These characteristic numbers (with complete proofs) were known classically.

\subsection{ Rational Plane Cubics}
Let $\M=\M_{0,0}(2,3)$.
$Pic(\M)\otimes \Q$ is freely generated by $\HH$ and the
unique boundary component $K$ corresponds to the degree partition
$1+2=3$. The algorithm described above yields the top intersections
of $\HH$ and $K$ inductively. Since $\tl{K}\subset \M_{0,1}(2,1)\times
\M_{0,1}(2,2)$, first the top intersections on these Kontsevich spaces
must be computed. $\M_{0,1}(2,1)$ has no boundary, hence all
top products are known. There is a unique boundary component
$B$ of $\M_{0,1}(2,2)$.
$\barr{B}\subset \M_{0,2}(2,1) \times \M_{0,1}(2,1)$.
Thus the top products on $\M_{0,2}(2,1)$ must be computed. Finally,
the unique boundary component of $\M_{0,2}(2,1)$ requires knowledge
of the top products on $\M_{0,3}(2,0)$ and $\M_{0,1}(2,1)$ which are
known. The answers for the top intersections of $\HH$ and $K$ on
$\M_{0,0}(2,3)$  ($dim \M_{0,0}(2,3)=8$) are:
$$\begin{array}{llr}
\HH^8 & & +12 \\
\HH^7  K & & +42 \\
\HH^6  K^2 & & +129 \\
\HH^5  K^3 & & + 285 \\
\HH^4  K^4 & & +336 \\
\HH^3  K^5 & & -(2541/ 4) \\
\HH^2  K^6 & & -(8259/ 16) \\
\HH    K^7 & &+ (19641/ 8) \\
K^8 & & - (44835/ 16)
\end{array}$$
Note since $K$ is $\Q$-Cartier, the intersections $\HH^i\cdot K^j$
need not be integers.
By Lemma (\ref{tan}), $\TT= {2\over 3} (\HH+K)$. The
characteristic number of plane cubics through $\alpha$ points
and tangent to $\beta$ lines is $\HH^\alpha \TT^\beta$:
$$\begin{array}{llr}
\HH^8 & & 12 \\
\HH^7\cdot \TT & & 36 \\
\HH^6 \cdot \TT^2 & & 100 \\
\HH^5 \cdot \TT^3 & & 240 \\
\HH^4\cdot \TT^4 & & 480 \\
\HH^3\cdot \TT^5 & & 712 \\
\HH^2\cdot \TT^6 & & 756 \\
\HH\cdot \TT^7 & & 600 \\
\TT^8 & & 400
\end{array}$$
These characteristic numbers have been calculated by
H. Zeuthen, S. Maillard, H. Schubert, G. Sacchiero,
S. Kleiman, S. Speiser, and P. Aluffi.
([S], [Sa], [K-S], [A]). Complete proofs appear in [Sa], [K-S], and [A].

\subsection{ Twisted Cubics in $\proj^3$}
In case $\M=\M_{0,0}(3,3)$, $Pic(\M)\otimes \Q$ is still generated
freely by $\HH$, $K$. A similar analysis yields the top intersections
($dim(\M)=12$):
$$\begin{array}{llr}
\HH^{12} & & +80160 \\
\HH^{11} K & & +121440 \\
\HH^{10} K^2 & & +148920 \\
\HH^9  K^3 & & +112080 \\
\HH^8 K^4 & & -7824 \\
\HH^7 K^5 & & -104100 \\
\HH^6 K^6 & & +35880 \\
\HH^5 K^7 & & + (190095/2) \\
\HH^4 K^8 & & - (222855/2) \\
\HH^3 K^9 & & -(674007/ 16) \\
\HH^2 K^{10} & & +(10112745/ 32) \\
\HH   K^{11} & & -(5995065/ 8) \\
K^{12} & & +(58086435/ 32)
\end{array}$$
The hyperplane tangency class is again $\TT={2\over 3}(\HH+K)$.
The characteristic number of twisted cubics through $\alpha$ lines
and tangent to $\beta$ planes is $\HH^{\alpha} \TT^{\beta}$:
$$\begin{array}{llr}
\HH^{12} & & 80160 \\
\HH^{11} \TT & & 134400 \\
\HH^{10} \TT^2 & & 209760 \\
\HH^9  \TT^3  & & 297280 \\
\HH^8 \TT^4 & & 375296 \\
\HH^7 \TT^5 & & 415360 \\
\HH^6 \TT^6 & & 401920 \\
\HH^5 \TT^7 & & 343360 \\
\HH^4 \TT^8 & & 264320 \\
\HH^3 \TT^9 & & 188256 \\
\HH^2 \TT^{10} & &128160 \\
\HH   \TT^{11} & & 85440 \\
\TT^{12} & & 56960
\end{array}$$
These characteristic numbers have been calculated by H. Schubert and others
([S], [K-S-X]). Complete proofs appear in [K-S-X].

\subsection{Rational Plane Quartics}
Let $\M=\M_{0,0}(2,4)$.
$Pic(\M)\otimes \Q$ is freely generated by
 $\HH$ and the
boundary components $J$, $K $ corresponding to the degree partitions
$2+2=4$, $1+3=4$. The top intersection numbers are
($dim\M_{0,0}(2,4)=11$):
$$\begin{array}{llrllllrllllr}
\HH^{11} &&  +620 \\
\HH^{10} K &&+1620 &&& \HH^{10}J &&+504 \\
\HH^9 K^2 &&+3564 &&& \HH^9JK  &&  +1512  &&&  \HH^9J^2   &&   +0\\
\HH^8 K^3 &&+4052 &&& \HH^8JK^2&&  +4536  &&& \HH^8J^2K &&   +0 \\
\HH^7 K^4 && -8340&&& \HH^7JK^3&&  +10920 &&& \HH^7J^2K^2 && +672 \\
\HH^6 K^5 && -48300 &&& \HH^6JK^4 && +15480 &&& \HH^6J^2K^3 && +4320 \\
\HH^5 K^6 && +1260 &&& \HH^5 J K^5 && -22296 &&& \HH^5J^2K^4 && +17184 \\
\HH^4K^7 && +153300 &&& \HH^4JK^6 && -22728  &&& \HH^4J^2K^5 && -11040 \\
\HH^3K^8 && -(338620/3)&&& \HH^3JK^7&&  +70056 &&& \HH^3J^2K^6 && -34560\\
\HH^2K^9&&  -(13690660/27)&&& \HH^2JK^8&&  +5880&&& \HH^2J^2K^7&&  +51072\\
\HH K^{10}&&   +(147582380/81)&&& \HH JK^9&& -385560&&&\HH J^2K^8 && +100800 \\
K^{11} &&  -(278947820/81)&&&  JK^{10} && +1310904 &&& J^2K^9 && -616896 \\
&&&&&&&&&&&&\\
\HH^8J^3&&      -364\\
\HH^7J^3K &&   -1260&&&  \HH^7J^4  &&    +630\\
\HH^6J^3K^2&&  -3852&&&  \HH^6J^4K &&  +1782&&&   \HH^6J^5 &&  -645\\
\HH^5J^3K^3&&  -8836&&&  \HH^5J^4K^2&&  +3588&&&  \HH^5J^5K &&  -(2385/2)\\
\HH^4J^3K^4&&  +4980&&&  \HH^4J^4K^3&&  -1788&&&  \HH^4J^5K^2&&  +906\\
\HH^3J^3K^5&&  +16356&&& \HH^3J^4K^4&&  -7830&&& \HH^3J^5K^3 &&  +(8241/2)\\
\HH^2J^3K^6&&  -22060&&& \HH^2J^4K^5&&  +7770&&& \HH^2J^5K^4 &&  -1815\\
\HH J^3K^7&&  -46452&&& \HH J^4K^6 &&   +22632&&& \HH J^5 K^5 &&   -(22125/2)\\
J^3K^8 && +255444&&& J^4K^7 &&-92232 &&& J^5K^6&&      +28920\\
&&&&&&&&&&&&\\
\HH^5J^6 &&     +(2419/8)\\
\HH^4J^6K &&   -(4743/8)&&&  \HH^4J^7 &&     +(765/2) \\
\HH^3J^6K^2&&-(18549/8)&&& \HH^3J^7K&& +1305&&&    \HH^3J^8&& -(5649/8)\\
\HH^2J^6K^3&&  -(3455/8)&&& \HH^2J^7K^2&&+(1923/2)&&&\HH^2J^8K&&-(6615/8)\\
\HH J^6K^4&& +(39075/8)&&& \HH J^7K^3&& -1680&&& \HH J^8K^2&&+(2163/8)\\
J^6K^5&&-(56631/8)&&& J^7K^4&& +(1701/2)&&& J^8K^3&&+(2289/8)\\
&&&&&&&&&&&&\\
\HH^2J^9&&  +(4375/8)\\
\HH J^9K&&  +189&&& \HH J^{10} && -(7875/32)\\
J^9K^2&&  -189&&&   J^{10}K && +0 &&& J^{11} && +(10143/128)
\end{array}$$
The line tangency class is $\TT={3\over 4}\HH+J+{3\over 4}K$.
The characteristic number of rational plane quartics through $\alpha$ points
and tangent to $\beta$ lines is $\HH^{\alpha} \TT^{\beta}$:
$$\begin{array}{llr}
\HH^{11} & & 620 \\
\HH^{10} \TT & & 2184 \\
\HH^{9}  \TT^2 & & 7200 \\
\HH^8  \TT^3  & & 21776 \\
\HH^7 \TT^4 & & 59424 \\
\HH^6 \TT^5 & & 143040 \\
\HH^5 \TT^6 & & 295544 \\
\HH^4 \TT^7 & & 505320 \\
\HH^3 \TT^8 & & 699216 \\
\HH^2 \TT^9 & & 783584 \\
\HH^1 \TT^{10} & &728160 \\
\TT^{11} & & 581904
\end{array}$$
\noindent The characteristic[ numbers of rational plane quartics
have been calculated by H. Zeuthen in [Z].

\subsection{Cuspidal Rational Plane Curves}
\label{cusp}
For $d\geq 1$, let $N_d$ be the number of irreducible,
nodal rational plane curves passing
through $3d-1$ general points in $\proj^2$. $N_d$ is a $\beta=0$
characteristic number. The numbers $N_d$ satisfy a beautiful recursion
relation ([K-M]):
$$N_1=1$$
$$\forall d\geq 2, \ \ N_d= \sum_{i+j=d,\  i,j>0}
N_i N_j i^2j \Bigg( j{3d-4\choose 3i-2} - i {3d-4\choose 3i-1}
\Bigg)\ \ \ .$$
The first few $N_d$'s are:
$$N_1=1, \ N_2=1, \ N_3=12, \ N_4=620, \ N_5=87304, \ N_6=26312976, \ \ldots$$
As a final application, the enumerative geometry of cuspidal rational
plane curves is considered. A rational plane curve, $C$, is {\em 1-cuspidal}
if the singularities of $C$ consist of nodes and exactly 1 cusp.
For $d\geq 3$, let $C_d$ be the
number of irreducible, 1-cuspidal rational plane curves passing
through $3d-2$ general points in $\proj^2$.
\begin{pr}
\label{cuspr}
The numbers $C_d$ can
be expressed in terms of the $N_d$:
$$\forall d\geq 3, \ \ C_d= {3d-3\over d} N_d \ + \ {1\over 2d}\cdot
\sum_{i=1}^{d-1} {3d-2\choose 3i-1}N_iN_{d-i} \big( 3i^2(d-i)^2 - 2 di(d-i)
\big) \ \ .$$
\end{pr}
\noindent The first few $C_d$'s are:
$$C_3=24, \ C_4=2304, \ C_5=435168, \ C_6=156153600, \ \ldots$$
$C_3$ is the degree of the locus of cuspidal cubics. $C_4$ has
been computed by H. Zeuthen ([Z]).
The 1-cuspidal numbers $C_d$ are evaluated by intersecting
divisors on $\M_{0,0}(2,d)$.

Let $d\geq 3$.
Let $M_{0,0}(2,d)$ be $\M_{0,0}(2,d)$ minus the boundary.
Let $Z \subset M_{0,0}(2,d)$ be the subvariety of maps that
are not immersions. It is easily seen $Z$ is of pure codimension
1 and the generic element of every component corresponds
to a 1-cuspidal rational plane curve. Let $\ZZ$ be the
Weil divisor obtained by the closure of $Z$ in $\M_{0,0}(2,d)$.
By the dimension reduction argument of section (\ref{cnum}),
the intersection cycle on $\M_{0,0}(2,d)$
\begin{equation}
\label{cy2}
\ZZ \cap \HH^{3d-2}
\end{equation}
determined by general points $P_1, \ldots, P_{3d-2}$
is of dimension (at most) 0 and lies in $Z$.
A simple modification of the corresponding argument
in section (\ref{cnum}) can be applied to show
(\ref{cy2}) is reduced for general choices of $P_j$.
Hence $C_d= \ZZ\cdot \HH^{3d-2}$.

The boundary of $\M_{0,0}(2,d)$ simply consists of
the $[{d\over2}]$ Weil divisors $K^i$ ($1\leq i \leq [{d\over 2}]$).
Recall $K^i$ is the boundary component corresponding to the degree
partition $i+(d-i)=d$. By Lemmas
(\ref{ihh}-\ref{i01}), the elements $\{\HH\} \cup \{K^i\}$
span a basis of $Pic(\M_{0,0}(2,d)) \otimes \Q$.
\begin{lm}
\label{zz}
The class of $\ZZ$ in $Pic(\M_{0,0}(2,d)) \otimes \Q$
is determined by ($d\geq 3$):
\begin{equation}
\label{laster}
\ZZ = {3d-3\over d}\HH \ + \ \sum_{i=1}^{[{d\over 2}]} {3i(d-i)-2d\over d} K^i
\ .
\end{equation}
\end{lm}

\begin{pf}
Let $S$, $\mu$, $\barr{S}$, $\lambda: C \rarr \M_{0,0}(2,d)$ be
exactly as in the proof of Lemma (\ref{tan}).
It will be checked that the left and right sides of
(\ref{laster}) have the same intersection with $C$.
As before, a point of the intersection $C\cdot K^i$ can arise
in two cases. A simple point of degree $i$ or $d-i$ can be blown-up.
Let $$C\cdot K^i= x_i + y_i$$
where $x_i$ and  $y_i$ are the number instances of the first and
second case respectively. Let $E_{x_i}$ be the union of the
$x_i$ exceptional divisors in $\barr{S}$ obtained from the
$x_i$ points of $C\cdot K^i$. Let $E_{y_i}$ be defined
similarly.

First, the intersection $C\cdot \ZZ$ is calculated.
Consider $\barr{\mu}:\barr{S} \rarr \proj^2$.
$\barr{\mu}^*(\oh_{\proj}(1))$ is the element
$(d,k) -\sum i E_{x_i} - \sum (d-i) E_{y_i}$.
The differential map yields an injection of {\em sheaves}:
\begin{equation}
\label{seq}
0 \ \rarr\  T_{\barr{S}}\  \stackrel{d\barr{\mu}}{\rarr}
 \barr{\mu}^*(T_{\proj^2}) \ \rarr\  Q\  \rarr \ 0.
\end{equation}
For general maps $\barr{\mu}$, $Q$ is line bundle supported
on a nonsingular curve $D$. The restriction of the
sequence (\ref{seq}) to $D$ yields an exact sequence of
{\em bundles} on $D$:
\begin{equation}
\label{seq2}
0 \ \rarr \ L \ \rarr\  T_{\barr{S}}|_D\  \stackrel{d\barr{\mu}}{\rarr}
 \barr{\mu}^*(T_{\proj^2})|_D \ \rarr\  Q|_D \  \rarr \ 0
\end{equation}
where $L$ is a line bundle on $D$.
Finally, there is exact sequence of bundles on
$D$ obtained from the projection $\barr{\pi}: \barr{S} \rarr C$:
\begin{equation}
\label{seq3}
0 \ \rarr\ V \ \rarr\  T_{\barr{S}}|_D \ \stackrel{d\barr{\pi}}{\rarr} \
\barr{\pi}^*(T_C) \rarr 0
\end{equation}
where $V$ is a line bundle on $D$. Maps in the
family $\barr{\pi}$ have zero differential
exactly at the points of intersection  $\proj(V) \cdot \proj(L) \subset
\proj(T_{\barr{S}}|D)$. Hence
$$C \cdot \ZZ = \proj(V) \cdot \proj(L)\ .$$
A lengthy, routine exercise in Chern classes and exact sequences now
yields:
$$C \cdot \ZZ = (6d-6)k + \sum_{i=1}^{[{d\over 2}]} (-3i^2+3i-2)x_i
 + \sum_{i=1}^{[{d\over 2}]} (-3(d-i)^2+3(d-i)-2)y_i.$$
Algebraic manipulation and the relation
$$C\cdot \HH = 2dk - \sum_{i=1}^{[{d\over 2}]} i^2 x_i -
\sum_{i=1}^{[{d\over 2}]} (d-i)^2 y_i$$
yields the result.
\end{pf}

It remains to compute $\ZZ\cdot \HH^{3d-2}$. By Lemma (\ref{zz}),
it suffices to determine the products $K^i \cdot \HH^{3d-2}$.
If $i\neq d/2$, the result
$$ K^i \cdot \HH^{3d-2} = {3d-2 \choose 3i-1} i(d-i)N_i N_{d-i}$$
is obtained from a simple geometric argument.
In case $i=d/2$,
division by 2 is required to account for symmetry:
$$ K^{d\over 2} \cdot \HH^{3d-2}= {1\over 2} {3d-2\choose 3{d\over 2}-1}
({d\over 2})^2 N^2_{d\over 2}\ .$$
Evaluation of $\ZZ \cdot \HH^{3d-2}$ yields the formula for $C_d$.
The proof of Proposition
(\ref{cuspr}) is complete.

\noindent Department of Math, University of Chicago, rahul@@math.uchicago.edu

\end{document}